\documentclass[12pt]{article}
\usepackage{amssymb}
\usepackage{amsmath}

\newcommand{\be}{\mbox{{\bf E}}}
\newcommand{\bp}{\mbox{{\bf P}}}
\newcommand{\eo}{E_{\omega}}
\newcommand{\hp}{\hat P}
\newcommand{\hq}{\hat Q}
\newcommand{\iot}{\int_{0}^{t}}
\newcommand{\ird}{\int_{\rd}}
\newcommand{\itq}{I_t^{(2)}(Q)}
\newcommand{\ot}{[0,t]}

\newcommand{\rd}{\R^d}

\newcommand{\tq}{\tilde Q}
\newcommand{\1}{{\bf 1}}

\newcommand{\wt}{W_{t}}

\newcommand{\zt}{Z_{t}}

\newcommand{\C}{\mathbb C}
\newcommand{\D}{\mathbb D}
\newcommand{\R}{\mathbb R}

\newcommand{\cf}{\mathcal F}
\newcommand{\ch}{\mathcal H}
\newcommand{\cl}{\mathcal L}

\newcommand{\cs}{\mathcal S}

\newcommand{\al}{\alpha}

\newcommand{\io}{\iota}
\newcommand{\la}{\lambda}
\newcommand{\om}{\omega}
\newcommand{\oom}{\Omega}

\newcommand{\si}{\sigma}

\newcommand{\vp}{\varphi}

\newcommand{\lp}{\left(}
\newcommand{\rp}{\right)}
\newcommand{\lc}{\left[}
\newcommand{\rc}{\right]}
\newcommand{\lcl}{\left\{}
\newcommand{\rcl}{\right\}}
\newcommand{\lln}{\left|}
\newcommand{\rrn}{\right|}
\newcommand{\lla}{\left\langle}
\newcommand{\rra}{\right\rangle}

\newcommand{\rat}{\right\rangle_{t}}

\newtheorem{theorem}{Theorem}[section]

\newtheorem{definition}[theorem]{Definition}
\newtheorem{example}[theorem]{Example}

\newtheorem{lemma}[theorem]{Lemma}

\newtheorem{proposition}[theorem]{Proposition}

\begin{document}

\thispagestyle{empty}
 \begin{center}
 {\Large\bf On the Brownian directed polymer}\\
{\Large\bf in a Gaussian random environment}

\vspace{0.5cm}

 by

\vspace{0.5cm}

{\bf Carles Rovira}\footnote{Partially done while the author was
visiting the Universit\'e Paris 13 with a CIRIT grant. Partially
supported by DGES grant BFM2003-01345 .\hfill}
and {\bf Samy Tindel}$^{2}$\\

\vspace{0.5cm}

 $^{1}$
{\it Facultat de Matem\`atiques,
Universitat de Barcelona,}
\\
\it Gran Via 585,
08007-Barcelona,
Spain
\\
{\it e-mail: carles.rovira@ub.edu}
\\
$^2$ {\it Institut \'Elie Cartan,
Universit\'e Henri Poincar\'e (Nancy),}
\\
{\it BP 239,
54506-Vandoeuvre-l\`es-Nancy,
France}
\\
{\it e-mail: tindel@iecn.u-nancy.fr}

 \end{center}

\vspace{1cm}

\begin{abstract}
In this paper, we introduce a model of Brownian polymer in a continuous
random environment. The asymptotic behavior of the partition function associated
 to this polymer measure is studied, and we are able to separate a weak and strong
 disorder regime under some reasonable assumptions on the spatial covariance of the
 environment. Some further developments, concerning some concentration inequalities
 for the partition function, are given for the weak disorder regime.
\end{abstract}

\vspace{1cm}


\vspace{1cm}


\newpage

\section{Introduction}
The directed polymer model in a random environment first appeared in the Mathematical Physics literature, as a canonical model of inhomogeneous systems (see e.g. \cite{DSp}, \cite{IS} for results in that direction). After some interesting relationships between this object and many other natural models of non-equilibrium dynamics have been established, the study of the polymer measure has been undertaken by Mathematicians, and a great amount of rigorous results is now available on the topic. These results concern basically the so-called partition function, the fluctuations and wandering exponents of the model, as well as the superdiffusive behavior of the polymer under the influence of the random media. On the other hand, a wide range of contexts have been explored: discrete random walks (see \cite{AZ}, \cite{Bo}, \cite{CH}, \cite{IS},
\cite{Si}), Brownian motion  in a discrete potential (see \cite{CO}, \cite{Co}), or Brownian motion in a Poisson-type potential (\cite{CY} or \cite{Wu}, \cite{Sn} for an undirected polymer).

\vspace{0.2cm}

This paper proposes to begin the study of a model which, from our point of view, is also worth considering, namely the Brownian polymer in a continuous Gaussian potential. More specifically, a complete description of our model can be given as follows:
\begin{enumerate}
\item
Our polymer will be modelized by a $d$-dimensional Brownian path
$\{\om_t;$ $t\ge 0  \}$, defined on a complete probability space
$(\hat\Omega,\hat\cf,\hat P)$ equipped with a filtration $\{\hat
\cf_t;t\ge 0  \}$. We will denote by $\eo$
 the expectation with respect to $\hat P$, that will be considered as the Wiener measure.
  We will also set $\hat P^x$ for the Wiener measure shifted by a constant $x\in\R^d$, which
  is of course the measure of a Wiener process with initial condition $x\in\R^d$.
\item
The random environment will be defined by a Gaussian landscape $B$ on $\R_+\times\R^d$,
 with rough fluctuations in time, and homogeneous with respect to the space coordinate: $B$
  will be given, on some probability space $(\Omega,\cf,{\bf P})$, as a centered Gaussian
   process whose covariance structure is
$$
\be\lc B(t,x) B(s,y) \rc=
\lp s\wedge t \rp Q(x-y),
$$
where $Q$ is a homogeneous covariance function such that $Q(0)<\infty$ (which
implies that $Q$ is bounded).
\end{enumerate}
Notice that $Q$ can also be represented by a Fourier transform
 procedure: there exists (see e.g. \cite{CV} for further details) a Gaussian
  independently scattered measure $M$ on $\R_+\times\R^d$ such that
$$
B(t,x)=\int_{\R_+\times\R^d}\1_{[0,t]}(s)e^{\io\la x}M(ds,d\la),
$$
where $\la x$ stands for the inner product of $\la$ and $x$ in $\rd$, and
where the law of $M$ is defined by the following covariance structure: for
any test functions $f,g:\R_+\times\R^d\to\C$, we have
\begin{multline*}
\be\lc \int_{\R_+\times\R^d} f(s,\la)M(ds,d\la)
\overline{\int_{\R_+\times\R^d}  g(s,\la) M(ds,d\la)} \,\, \rc\\
= \int_{\R_+\times\R^d} f(s,\la) \overline{g(s,\la)}\hat Q(d\la)
ds,
\end{multline*}
and the finite (real) measure $\hq$ is the Fourier transform of $Q$. With
this notation in mind, we can complete the description of our polymer measure by
\begin{itemize}
\item[3.]
For any $t>0$, the energy of a given path (or configuration) $\om$ on $[0,t]$ will be given by
$$
-H_t(\om)=\iot B(ds,\om_s)
=\iot\int_{\R^d}e^{\io\la \om_s}M(ds,d\la).
$$
Notice that, for any fixed path $\om$, $H_t(\om)$ is a centered Gaussian random variable with
variance $tQ(0).$

Based on this Hamiltonian, for any $x\in\R^d$, and a given
constant $\beta$ (interpreted as the inverse of the temperature of
the system), we will define our (random) polymer measure by
$$
dG_t^x(\om)=\frac{e^{-\beta H_t(\om)}}{Z_t^x}d\hat P^x(\om) ,\quad
\mbox{ with }\quad Z_t^x=\eo^x\lc e^{-\beta H_t(\om)} \rc.
$$
\end{itemize}
In the sequel, we will also consider the Gibbs average with respect
to the polymer measure, defined as follows: for all  $t\ge 0$, $n\ge 1$,
and for any bounded measurable functional $f:(C(\ot;\rd))^n\to\R$, we set
\begin{equation}\label{avpol}
\lla f\rra_t=\frac{\eo^x\lc f(\om^1,\ldots,\om^n)
e^{-\beta\sum_{l\le n}H_t(\om^l)}\rc}{Z_t^n},
\end{equation}
where the $\om^l, 1 \le l \le n,$ are understood as independent
Brownian configurations.
\vspace{0.2cm}

Our paper will be mainly concerned with the study of the partition function $\zt$ of the model described above, and let us mention already that, for the results we have obtained so far, the relevant parameters for our model will be the covariance function $Q$, and the inverse of the temperature $\beta$. Based on these parameters, we will get the following results:
\begin{itemize}
\item
A concentration inequality and the almost sure limit for $\frac{1}{t}\log(\zt)$.
\item
A natural definition of the weak and strong disorder regime
for our polymer (see Definition \ref{defweak})
\item
In the case of a covariance function $Q(x)$ that
can be written as $\tilde Q(|x|)$ with $\tilde Q:\R_+\to\R_+$, for $d\ge 3$,
we will show that a sufficient condition in order to be
in the weak disorder regime is
$\int_0^\infty u\tilde Q(u)du<\infty$ and $\beta$ small enough.
\item
For the weak disorder regime, we will show some refined concentration results for $\frac{1}{t}\log(\zt)$, using some general techniques taken from the random media literature (cf \cite{Tbk}, \cite{CH}).
\item
We will show that, for
 any $d\ge 1$, if $c_1 (1 + |x|^{2})^{-\lambda} \le Q(x) \le c_2
 (1+ |x|^2)^{-\hat\lambda}$ for some constants
 $c_1>0, c_2 >0$ and $0 < \hat\lambda \le \lambda < \frac12$, then the polymer will be in the strong disorder
 regime, regardless of the value of $\beta>0$.
\end{itemize}
Of course, many problems remain open for this model:
behavior of the wandering and fluctuation exponents, existence of a
covariance function $Q$ for which a phase transition can be seen as $\beta$
grows to $\infty$, computations involving the overlap function associated to the
 model (which will be defined by equation (\ref{defoverlap})), etc. We plan to report on these
  issues in a subsequent paper.

\vspace{0.2cm}

It is also worth mentioning that we have chosen to deal with this specific model for two main reasons:
\begin{enumerate}
\item
The continuous Gaussian model, which is physically a rea\-so\-na\-ble choi\-ce, allows us to use the huge amount of techniques available for this kind of processes (stochastic calculus, concentration inequalities, Malliavin calculus, among others), leading to some quite simple proofs of the main results contained in this paper.
\item
It is well known that $\zt$ behaves, in law, like the Feynman-Kac representation of $u(t,0)$, where $u(t,x)$ is the mild solution to the stochastic PDE
\begin{equation}\label{heateq}
\partial_t u(t,x)=\Delta u(t,x)+\beta u(t,x) \dot W(dt,dx),\quad
t\ge 0, x\in\R^d,
\end{equation}
understood in the Stratonovich sense, with $u(0,x)=1$, and thus
$$
\lim_{t\to\infty}\frac{1}{t}\log(\zt)
$$
can be interpreted as the Lyapounov exponent for this equation. Our problem is thus closely related to the one considered in \cite{CV}, \cite{TV1}, \cite{TV2} (see  also \cite{CM} and \cite{CS} for the discrete case), and though the questions adressed here are not exactly the same as in the latter papers, we believe that the present article gives some more insight on the (rather) old problem of the Lyapounov exponent for equation (\ref{heateq}). For instance, to our knowledge, the existence of this Lyapounov exponent had never been proven before, and its exact computation for $d\ge 3$ had never been performed either.
\end{enumerate}

Our paper will be organized as follows: at Section 2, we recall some basic notions and theorems of stochastc analysis that will be used in the sequel, we give some results on the alomost sure behavior of $\frac{1}{t}\log(\zt)$, and we define our notions of weak and strong disorder. At Section 3, we study in detail the weak disorder regime. At Section 4, we give a basic example of a strong disorder situation.

\section{Almost sure limit of the partition function}
In this section, we will give some basic results about the almost sure
 convergence of $Z_t$, and some rough bounds  on its limit. This will allow
 us to define precisely a notion of weak and strong disorder for the polymer
 measure. First of all, we will introduce some notation on Malliavin calculus
  for the Gaussian measure $M$, that we will use throughout the paper.

\subsection{Malliavin calculus preliminaries}

We will give here some notations and basic results, taken mainly
from \cite{Mal}, \cite{Nu} and \cite{U}.  Let us specify first the
Wiener space we will consider: for any test functions
$f,g:\rd\to\C$, set
$$
(f,g)=\int_{\rd\times\rd}f(\la) \overline{g(\la)} \hq(d\la).
$$
Call $H$ the completion of $C_c^\infty(\rd)$ with respect to that
positive bilinear form, and $(\cdot,\cdot)_H$ the corresponding
inner product. Set also $\ch=L^2(\R_+;H)$. The Gaussian process
$M$ can be seen as a zero-mean Gaussian family $\{ M(h); \, h\in
\ch \}$ satisfying
$$
\be\lc M(h_1)\overline{M(h_2)} \rc = (h_1,h_2)_{\ch} \equiv
\int_{\R_+}\lp h_1(t), h_2(t) \rp_H dt, \qquad h_1,h_2\in \ch,
$$
where we have set, for $h\in \ch$,
$$
M(h)=\int_{\R_+\times\rd} h(\si,\la) M(d\si,d\la).
$$
Furthermore, we will assume that $\cf$ is generated by $M$. Then $(M,\ch,\bp)$
defines a Wiener space on $C(\R_+\times\rd;\C)$, on which the traditional tools
 of Malliavin calculus can be introduced. Let us recall some of them for sake of
  completeness: a smooth functional of $M$ will be of the form
\begin{equation}\label{smoothfct}
F=f\lp M(u_1),\ldots,M(u_m) \rp,\qquad
m\ge 1, u_j\in\ch, f\in C^\infty(\C^m),
\end{equation}
and we will denote by $\cs$ the set of such functionals. Now, for $F$
 as in (\ref{smoothfct}), the Malliavin derivative of $F$ will be defined,
  as an element of $\ch$, by
\begin{equation}\label{dfbsmooth}
D_{t,\la}F=
\sum_{j=1}^{m}\partial_{x_j}f\lp M(u_1),\ldots,M(u_m) \rp
u_j(t,\la).
\end{equation}
Then it can be shown that the operator $D:\cs\to\ch$ is closable, and, as usual,
for any $p\in[1,\infty]$, we will denote by $\D^{1,p}$ the Sobolev space obtained by
completing $\cs$ with respect to the norm
$$
\| F\|_{1,p}^{p}=
\be\lc |F|^p \rc + \be\lc |DF|_{\ch}^p \rc.
$$
Notice that the following chain rule is available for functionals
$F$ in $\D^{1,p}$: if $\vp:\R\to\R$ is a smooth function such
 that $\vp(F),\vp'(F)\in L^{\al}(\oom)$ for any $\al>0$,
 then $\vp(F)\in \D^{1,r}$ for any $r<p$, and
\begin{equation}\label{chrule}
D\vp(F)=\vp'(F)\, DF.
\end{equation}
Let us also mention that, among all the elaborated
 integration by parts formulae of the Malliavin calculus, we will
  only use the following basic one in the sequel: if $F\in\D^{1,2}$ and
   $u$ is a deterministic element of $\ch$, then
\begin{equation}\label{tibp}
\be\lc F\, M(u)\rc=
\be\lc \lp DF,u \rp_\ch \rc.
\end{equation}

Concentration inequalities are a useful tool in random system theory, and
 we will use the following one, taken from \cite{U}:
\begin{proposition}\label{concentracio}
Let $F\in\D^{1,p}$ for some $p>1$, and suppose that $DF\in L^\infty(\oom;\ch)$. Set
 $m=\be[F]$ and $\si^2=\|DF\|_{\infty}$. Then we have, for any $c
 \ge 0$,
$$
\bp\lp |F|\ge c \rp
\le
2\exp\lp -\frac{(c-m)^2}{2\si^2} \rp.
$$
\end{proposition}

\vspace{0.5cm}

\noindent
We will also need a refinement of Proposition \ref{concentracio}, for which
 we have to introduce a little more notation and a  0-1 type law, that we learned
 from \cite{U2}:

 \begin{lemma}\label{lemanou} Let $J$ be a measurable set in $M$, such that $J +
 \ch \subset J$. Then $P(J)\in \{0,1\}.$
\end{lemma}

\vspace{0.3cm}

\noindent {\bf Proof:} It is well known (see \cite[page 31]{Nu})
that $F=\1_J$ is an element of $\D^{1,2}$ iff $P(J) \in \{0,1 \}$.
However if $J + \ch \subset J$, we also have $J + \ch = J$, and
this easily yields $D \1_J = 0$ and thus $\1_J \in \D^{1,2}.$
 \hfill $\Box$

\bigskip

For $h\in\ch$, set
  $\tilde h_t=\iot h_s ds$. For a measurable
   subset $A$ of $(M,\ch,\bp)$ and $m\in M$, define
$$
q_A(m)=\inf\lcl |h|_{\ch};\, m+\tilde h\in A \rcl.
$$
Then we claim that:
\begin{lemma}\label{ltala}
Suppose $\bp ( A ) \ge p > 0$. Then, for any  $u>0$,
$$
\bp \Big( q_A  > u + (2 \log(2/p))^{\frac12} \Big) \le 2 \exp
\Big( -\frac{u^2}{2} \Big).
$$
\end{lemma}

\vspace{0.4cm}

\noindent {\bf Proof:} Let us prove first that $|Dq_A|_{\ch}\le 1$
(see \cite{FU} for further details on this functional). First, by
Lemma \ref{lemanou}, $q_A$ is almost surely finite. Indeed, if
$J=\{\om; q_A(\om)<\infty  \}$, then it is easily checked that $J
+ \ch \subset J$, and thus $\bp(J)\in\{0,1 \}$. On the other hand,
$A\subset J$, and by assumption $\bp(A)>0$, which gives the
finiteness of $q_A$. Furthermore, if $l\in\ch$,  by the usual
triangular inequality, we have almost surely
$$
q_A(m+\tilde l)\le q_A(m)+|l|_\ch,
$$
and hence $q_A$ is a Lipschitz functional on $M$ with Lipschitz
constant 1, which yields, in particular, $|Dq_A|_{\ch}\le 1$. Then
applying Proposition \ref{concentracio}, we obtain that for any $u
>0,$
$$
\bp \left( \left\vert q_A - \be [q_A] \right\vert > u \right) \le 2 \exp
\left( - \frac{u^2}{2} \right).
$$
The proof now follows the lines of Talagrand \cite{Tbk} (see also \cite{CH}):
if $u <\be [q_A],$ we have
$$
p \le  \bp(A) \le \bp \left( \left\vert q_A - \be [q_A] \right\vert
> u \right) \le 2 \exp \left( - \frac{u^2}{2} \right),
$$
which yields that $\be [q_A] \le (2 \log(2/p) )^\frac12$, and the
proof is now easily completed.

\hfill $\Box$

\subsection{Almost sure behavior}
Let us begin with a Markov type decomposition for $Z_t^x$: for
$x,y\in\rd$, and $t,h\ge 0$, set
\begin{eqnarray*}
Z_{t,t+h}(x,y,B) &=&E_{\om}^{x}\left[
e^{\beta\int_{0}^{t+h}B(ds,\omega
_{s})}\left|\, \omega _{t}=y \right.\right] , \\
Z_{t}(x,B) &=&E_{\om}^{x}\left[ e^{\beta\int_{0}^{t}B(ds,\omega _{s})}%
\right] .
\end{eqnarray*}
Then the following property holds true:
\begin{lemma}\label{ch1}
Let $p_t$ be the heat kernel on $\rd$ at time $t\ge 0$, and set, for $t,s\in\R_+$ and
 $x\in\rd$,
$\theta _{t}B(s,x)=B(s+t,x).$
Then, for any $t,h\ge 0$, we have
$$
Z_{t+h}(x,B) = \int_{\R^{d}}Z_{h}(y,\theta
_{t}B)Z_{t,t}(x,y,B)p_{t}(dy).
$$
\end{lemma}

\vspace{0.5cm}

\noindent {\bf Proof:} Notice first the relationship
\begin{eqnarray*}
Z_{t+h}(x,B) &=&\int_{\rd }Z_{t,t+h}(x,y,B)p_{t}(dy) \\
&=&\int_{\rd }E_{\om}^{x}\left[ e^{\beta\int_{0}^{t}B(ds,\omega
_{s})}e^{\beta\int_{t}^{t+h}B(ds,\omega _{s})}\Big| \omega
_{t}=y \right] p_{t}(dy).
\end{eqnarray*}
Thus, we have that
\begin{eqnarray*}
Z_{t+h}(x,B) &=&E_{\om}^{x}\left[ E_{\om}^{x}\left[
e^{\beta\int_{0}^{t}B(ds,\omega
_{s})}e^{\beta\int_{t}^{t+h}B(ds,\omega _{s})} \Big| \hat \cf_{t}
  \right] \right]\\
&=&E_{\om}^{x}\left[  e^{\beta\int_{0}^{t}B(ds,\omega
_{s})}E^{\omega
_{t}}_{\widehat{\omega }}\left[ e^{\beta\int_{0}^{h}(\theta _{t}B)(ds,\widehat{%
\omega }_{s})}\right]  \right] \\
&=&\int_{\rd }E_{\om}^{x}\left[ e^{\beta\int_{0}^{t}B(ds,\omega _{s})}E^{y}_{%
\widehat{\omega }}\left[ e^{\beta\int_{0}^{h}(\theta _{t}B)(ds,\widehat{\omega }%
_{s})}\right] \Big| \omega _{t}=y  \right] p_{t}(dy) \\
&=&\int_{\rd }Z_{h}(y,\theta _{t}B)Z_{t,t}(x,y,B)p_{t}(dy),
\end{eqnarray*}%
where  $\hat \omega$ denotes a $d$-dimensional Brownian path, independent of $\om$.

\hfill $\Box$

\vspace{0.5cm}

As usual in disordered systems theory, the free energy, defined by
$$
p_t (\beta) = \frac{1}{t} \be  \left[ \log \lp Z_t^x \rp \right],
$$
will play an important role in the qualitative description of the asymptotic behavior of the system.
Note that, taking into account the space homogeneity of $B$, $\be[\log(Z_t^x)]$ will be independent of the parameter $x\in\rd$. This is why we will concentrate now on this quantity for $x=0$. In fact, from now, $x$ will be understood as 0 when not specified, and $\eo,\zt$ will stand for $\eo^0,\zt^0$, etc.
We are now in position to state a first basic result about the limit of the quantity $p_t (\beta)$.
\begin{proposition}\label{p14}
For all $\beta >0$ there exists a constant $p(\beta)>0$ such that
\begin{equation*}
p(\beta)\equiv\lim_{t \to \infty} p_t (\beta) = \sup_{t \ge 0}
p_t(\beta) .
\end{equation*}
\end{proposition}

\vspace{0.5cm}

\noindent {\bf Proof:}
This result is presumably fairly standard, but we include its proof for sake of readability: for $t,h\ge 0$, invoking Lemma \ref{ch1}, Jensen's inequality, and the independence of the time increments of $B$, we get
\begin{align*}
&\be  \left[ \log Z_{t+h} (0,B) \right]\\
& =
\be  \left[ \log \int_{\R^d} Z_{h} (y, \theta_t B) Z_{t,t} (0,y,B) p_t(dy) \right] \\
& =  \be  \left[ \log Z_{t} (0,B) \right]
+ \be  \left[ \log \int_{\R^d} Z_{h} (y, \theta_t B) \frac{Z_{t,t} (0,y,B)}{Z_t(0,B)} p_t(dy) \right] \\
& \ge
 \be  \left[ \log Z_{t} (0,B) \right]
+ \be  \left[  \int_{\R^d}  \big( \log Z_{h} (y, \theta_t B) \big) \frac{Z_{t,t} (0,y,B)}{Z_t(0,B)} p_t(dy) \right] \\
& =   \be  \left[ \log Z_{t} (0,B) \right] + \int_{\R^d}    \be  \left[
\log Z_{h} (y, \theta_t B)  \right]
 \be  \left[  \frac{Z_{t,t} (0,y,B)}{Z_t(0,B)} \right] p_t(dy) \\
& =   \be  \left[ \log Z_{t} (0,B) \right] +   \be  \left[  \log Z_{h} (0,
\theta_t B)  \right]
 \be  \left[  \int_{\R^d}  \frac{Z_{t,t} (0,y,B)}{Z_t(0,B)}  p_t(dy)  \right] \\
& =   \be  \left[ \log Z_{t} (0,B) \right] +   \be  \left[  \log Z_{h} (0,
 B)  \right].
\end{align*}
Notice that, in the above inequality, we have also used the fact that,
for any $y\in\rd$,
$$
\be[\log Z_{h} (y,\theta_t B)]=\be[\log Z_{h} (0,\theta_t B)],
$$
thanks to the space homogeneity of $B$.
Thus, for all $t,h\ge 0$,
$$
(t+h)p_{t+h}(\beta)\ge t p_{t}(\beta)+ h p_{h}(\beta).
$$
This easily yields
$$
\lim_{t \to \infty} p_t (\beta) = \sup_{t \ge 0}
p_t(\beta):=p(\beta)
$$
by a superadditivity argument.

\hfill
$\Box$

\vspace{0.5cm}

We will now summarize some elementary properties of $p(\beta)$:
\begin{proposition}
The function $p$ introduced at Proposition \ref{p14} satisfies:
\begin{enumerate}
\item
The map $\beta\mapsto p(\beta)$ is a convex nondecreasing function
on $\R_+$.
\item
The following upper bound holds true:
\begin{equation}\label{roughbnd}
p(\beta)\le \frac{\beta^2}{2} Q(0).
\end{equation}
\item
$\bp$-almost surely, we have
\begin{equation}\label{aslimfree}
\lim_{t \to \infty} \frac{1}{t} \log Z_t  = p(\beta).
\end{equation}
\end{enumerate}
\end{proposition}

\vspace{0.5cm}

\noindent
{\bf Proof:}
We will divide this proof in several steps.

\vspace{0.5cm}

\noindent
{\bf Step 1:}
The convexity of $p$ is a trivial consequence of H\"older's inequality.

\vspace{0.5cm}

\noindent
{\bf Step 2:}
In order to prove the third point of the proposition, let us compute
 first the Malliavin derivative of $U_t\equiv\frac{1}{t}\log(Z_t)$: since
  $\zt$ is $\cf_t$-adapted, we have $D_{\tau,\la}U_t=0$ if $\tau> t$, and
  if $\tau\le t$, according to (\ref{chrule}), we get $U_t\in\D^{1,2}$ and
$$
\frac{D_{\tau,\la}(\log(\zt))}{t}
=
\frac{D_{\tau,\la}(\zt)}{t \zt}
=\frac{\beta \eo\lc e^{\io\la\om_{\tau}}e^{-\beta H_t(\om)} \rc}{t\zt}
=\frac{\beta\lla e^{\io\la\om_{\tau}}\rat}{t},
$$
where we have used (\ref{dfbsmooth}) in order to differentiate $e^{-\beta H_t(\om)}$.
Notice also that the order of $D$ and $\eo$ can be interchanged by a simple uniform
convergence argument. Hence, by definition of the inner product in $\ch$,
\begin{eqnarray*}
|DU_t|_{\ch}^{2}&=&
\lp\frac{\beta}{t}\rp^2
\iot\int_{\rd} \lla e^{\io\la\om_{\tau}}\rat
\lla e^{-\io\la\om_{\tau}}\rat \hq(d\la) d\tau\\
&=&
\lp\frac{\beta}{t}\rp^2
\iot\lla \int_{\rd} e^{\io\la(\om_{\tau}^1-\om_{\tau}^2)}\hq(d\la)\rat
 d\tau\\
&=&
\lp\frac{\beta}{t}\rp^2 \iot\lla Q\lp\om_{\tau}^1-\om_{\tau}^2\rp\rat
 d\tau.
\end{eqnarray*}
Observe that, in the above expression, $\om^1$ and $\om^2$ are
understood as two independent configurations under the polymer
measure, and that we have also used the notation (\ref{avpol}). In
particular,
$$
|DU_t|_{\ch}^2\le \frac{\beta^2Q(0)}{t},
$$
almost surely. Hence, as a direct consequence of Proposition
\ref{concentracio}, we get
\begin{equation}\label{concentzt}
P \left( \left| \frac{1}{t} \log Z_t  - p_t (\beta) \right|
>c \right)
\le
2 \exp \left( - \frac{t c^2}{4Q(0)\beta^2} \right),
\end{equation}
from which (\ref{aslimfree}) can be deduced by a standard Borel-Cantelli argument.

\vspace{0.5cm}

\noindent
{\bf Step 3:}
In order to prove the bound (\ref{roughbnd}), let us just observe that
Jensen's inequality trivially yields
$$
p_t (\beta) \le \frac{1}{t} \log \be [Z_t  ] .
$$
However, the computation of $\be [Z_t (\beta) ]$ is an easy task: for any
fixed $\om$, $-H_t(\om)$ is a Gaussian random variable, and hence
\begin{equation}\label{espht}
\be\lc e^{-\beta H_t(\om)} \rc= \exp\lp \frac{\beta^2
\be[(H_t(\om))^2]}{2}\rp =e^{\frac{\beta^2 Q(0)t}{2}}.
\end{equation}
Thus, for any $\beta\ge 0$,
\begin{equation}
p_t (\beta) \le \frac{\beta^2 Q(0)}{2}.
\end{equation}

\hfill $\Box$

\subsection{Weak and strong disorder}

The amount of influence of the environment $B$ on the
 path $\om$ is usually captured through the behavior of $\zt$
 (see \cite{CY}, \cite{CH}). More specifically, we can argue as
 follows: recall that relation (\ref{roughbnd}) states
  that $p(\beta)\le \frac{\beta^2Q(0)}{2}$. The weak disorder regime
  is then naturally characterized by the relation
$$
p(\beta)= \frac{\beta^2Q(0)}{2}
\quad \mbox{ i.e. }\quad
\lim_{t\to\infty}\frac{1}{t}\be\lc \log(\zt) \rc
=
\lim_{t\to\infty}\frac{1}{t}\log\lp\be\lc \zt \rc \rp,
$$
while the strong disorder phase should be defined by
$p(\beta)< \frac{\beta^2Q(0)}{2}$. However, it will be more convenient
to define the weak and strong disorder regimes through an associated process:
set, for $t\ge 0$,
\begin{equation}\label{defwt}
W_{t}=Z_{t} \exp\lp-\frac{\beta^{2}Q(0) t }{2}\rp.
\end{equation}
Then it is easily seen that $W$ is a positive $\cf_t$-martingale,
that converges almost surely. Set then
$$
W_\infty=\lim_{t\to\infty}W_t.
$$
By Kolmogorov's 0-1 law and an easy elaboration of \cite[Lemma
2]{Bo} , we have
$$
\bp\lp W_\infty=0 \rp \in\{ 0,1 \}.
$$
Observe that if $W_\infty>0$ almost surely, then $\log(W_\infty)$
is finite almost surely, and hence
$$
\mbox{a.s.}-\lim_{t\to\infty} \Big( \log(\zt)
-\frac{\beta^{2}tQ(0)}{2} \Big) =\log(W_\infty),
$$
which yields
$$
\mbox{a.s.}-\lim_{t\to\infty}\frac{\log(\zt)}{t}
=\frac{\beta^{2}Q(0)}{2},
$$
and hence
\begin{equation}\label{winfpbet}
p(\beta)=\frac{\beta^{2}Q(0)}{2}.
\end{equation}
In other words, $W_\infty>0$ implies $p(\beta)=\frac{\beta^{2}
Q(0)}{2}$, and hence a weak disorder type behavior of the polymer.
This is why we will adopt the following definition:
\begin{definition}\label{defweak}
We will say that the polymer is in a strong disorder regime if $W_\infty=0$ almost
surely, while the weak disorder phase will be defined by $W_\infty>0$ almost surely.
\end{definition}

Another relevant quantity for the study of disordered systems is the so-called
overlap, that measures the similarity of two independent configurations under
 the considered random measure. In our case, this overlap is of the form
\begin{equation}\label{defoverlap}
\frac{1}{t}\iot\lla Q(\om^1_s-\om^2_s) \rra_s ds,
\end{equation}
and observe that, since $Q(x)$ is usually a decreasing function of $|x|$, the last quantity really measures how close $\om^1$ is from $\om^2$. One is then also allowed to relate the behavior of $\wt$ and of the overlap in the following way:
\begin{proposition}\label{p51}
Let $\wt$ be defined by (\ref{defwt}) for $t>0$, and consider the statements:
\begin{enumerate}

\item $W_\infty > 0$ almost surely.

\item $\int_0^\infty \langle Q(\omega_s^1 - \omega_s^2) \rangle_s
ds < \infty.$

\item $L^1 - \lim_{t \to \infty} W_t = W_\infty.$

\end{enumerate}

Then 1. and 2. are equivalent, and are both a
consequence of statement 3.

\end{proposition}

\vspace{0.5cm}

\noindent {\bf Proof:}
Let us check first that 3. implies 1. The
convergence in $L^1$ that we are assuming implies that
$$
E( W_\infty ) = \lim_{t \to \infty} E(W_t) =1.
$$
Using Kolmogorov's 0-1 law, we get $P(W_\infty > 0)=1$.

\vspace{0.5cm}

\noindent
Let us prove now the equivalence between 1 and 2: for $t\ge 0$, set
$$
N_t=\beta\iot\ird e^{\iota\la\om_s}M(ds,d\la).
$$
Then, for any fixed configuration $\om$, $N_t$ is a martingale,
whose quadratic variation process is given by
$$
[N]_t=\beta^2\iot\ird \hq(d\la)=\beta^2 Q(0)t.
$$
Furthermore, we have
$$
\wt=\eo\lc \exp\lp N_t-\frac{\beta^2 Q(0)t}{2}\rp \rc,
$$
and It\^o's formula applied to $\vp(x)=e^x$ gives the following martingale
 decomposition for $W$:
\begin{equation}\label{mgw}
\wt=1+\beta
\iot\ird\eo
\lc  e^{\iota\la\om_s} \exp\lp N_s-\frac{\beta^2 Q(0)s}{2}\rp\rc
M(ds,d\la),\quad t\ge 0.
\end{equation}
The process $\wt$ is also almost surely strictly
positive. Hence, one can apply again It\^o's formula
to the function $\psi(x)=\log(x)$ to get
$$
\log(W_t) =  \int_0^t \frac{dW_s}{W_s} - \frac12 \int_0^t
\frac{d[W]_s}{W_s^2} = M_t- \frac12 A_t,
$$
with
\begin{eqnarray*}
M_t&=&
\beta\iot\ird\lla e^{\iota\la\om_s} \rra_s M(ds,d\la)\\
A_t&=& \beta^2 \iot\ird\lla e^{\iota\la(\om^1_s-\om^2_s)} \rra_s
\hq(d\la) ds =\beta^2 \iot \lla Q\lp \om^1_s-\om^2_s \rp \rra_s
ds.
\end{eqnarray*}
Moreover, notice that
$\{M_t;t \ge 0\}$ is a martingale with quadratic variation $A_t$,
and that we can write
\begin{equation}\label{desco}
\log(W_t) = A_t \left( \frac{M_t}{A_t} - \frac12 \right).
\end{equation}
Now one can argue as follows:
\begin{description}

\item (a) Assume that $A_\infty = \infty$.
We can now apply the strong law of large numbers for continuous
martingales (see for instance \cite{RY}),  that
implies $\frac{M_t}{A_t} \to 0$,
almost surely. Then expression (\ref{desco}) gives us that $W_\infty=0$.

\item (b) Assume that $A_\infty < \infty$. If $A_\infty < \infty$ then we have that $M_t$ is
a $L^2$-bounded martingale which converges almost surely to $M_\infty$ when
$t$ goes to $\infty$. So $M_\infty < \infty$ almost surely, which clearly
yields that $\log W_\infty > - \infty$ almost surely.

\end{description}

\hfill $\Box$

\section{The weak disorder regime}

In this section, we will give a sufficient condition
under which the polymer is in the weak disorder phase. It
 is usually satisfied when $\om$ is transient under $\hat P$,
 and $Q$ rapidly decays at infinity, as we will see in some examples.
 Eventually, we will show that the
  concentration of $\frac{1}{t}\log(\zt)$ below $\frac{\beta Q(0)}{2}$ occurs at a
  higher speed than the one indicated by
(\ref{concentzt}). For sake of readability, we will make, in this
section,
 the following assumption:
\begin{itemize}
\item[{\bf (H)}]
$Q$ is a symmetric function from $\rd$ to $\R$ and $\beta$ a positive constant satisfying
$$
E_\omega  \left[ e^{\frac{\beta^2}{2} I_{\infty}(Q)} \right] <
\infty, \quad \mbox{ where } \quad I_{\infty}(Q)=\int_0^\infty Q(
\omega_s) ds
$$
\end{itemize}
Let us begin with our sufficient condition ensuring the weak disorder regime:
\begin{proposition}\label{pp}
Under hypothesis (H), we have
$$
\bp(W_\infty >0)=1 \quad\mbox{ and } \quad
p(\beta) = \frac{\beta^2 Q(0)}{2}  .
$$
\end{proposition}

\vspace{0.5cm}

\noindent {\bf Proof:}
We will divide this proof in two steps

\vspace{0.5cm}

\noindent {\bf Step 1:}
Let us compute $\be[\zt^2]$: notice that
$$
Z_{t}^{2}=E_{\omega }\left[ e^{\beta \int_{0}^{t}B(ds,\omega
_{s}^{1})+B(ds,\omega _{s}^{2})}\right],
$$
and hence, using the fact that
$\int_{0}^{t}B(ds,\omega_{s}^{1})+B(ds,\omega _{s}^{2})$ is a
Gaussian random variable for any fixed $\om^1$ and $\om^2$, we get
$$
\be \left[Z_{t}^{2}\right] = E_{\omega }\left[ \exp\lp
\frac{\beta ^{2}}{2}\be\left[ \left| \int_{0}^{t}B(ds,\omega
_{s}^{1})+B(ds,\omega _{s}^{2})\right| ^{2}\right] \rp\right].
$$
On the other hand, since
$$
\iot B(ds,\omega _{s}^{1})+B(ds,\omega _{s}^{2})= \iot\ird\lp
e^{\iota \la\om_s^1}+e^{\iota \la\om_s^2} \rp M(ds,d\la) ,
$$
we get
\begin{align*}
&\be\left[ \left| \int_{0}^{t}B(ds,\omega
_{s}^{1})+B(ds,\omega _{s}^{2})\right| ^{2}\right]\\
&=
\iot\ird \lp e^{\iota \la\om_s^1}+e^{\iota \la\om_s^2} \rp
\lp e^{-\iota \la\om_s^1}+e^{-\iota \la\om_s^2} \rp
\hq(d\la) ds\\
&=2\lp Q(0)t+\iot Q(\om_s^1-\om_s^2) ds\rp.
\end{align*}
Thus
\begin{equation}\label{sqzt}
\be \left[ Z_{t}^{2}\right] =E_{\omega }\left[ e^{\beta ^{2} ( Q(0)
t +  \int_{0}^{t}Q(\omega _{s}^{1}-\omega _{s}^{2})ds )} \right] .
\end{equation}

\vspace{0.5cm}

\noindent {\bf Step 2:}
Recall now that $\wt=\zt e^{2^{-1}\beta Q(0) t}$. Thus, using the fact
that $\om^1-\om^2$ can be written, in law, as $2^{1/2}\om$, where $\om$ is
 again a $\hp$-Brownian motion, we get
$$
\be \left[ W_{t}^{2}\right] = E_{\omega }\left[ e^{\beta ^{2}
\int_{0}^{t}Q(\omega _{s}^{1}-\omega _{s}^{2})ds} \right] \le
E_{\omega }\left[ e^{\frac{\beta^2}{2} I_{\infty}(Q)} \right].
$$
Hence, under assumption (H), $\wt$ is a bounded martingale in
$L^2$ with $\be [W_t]=1$, which yields in particular
$\be[W_\infty]=1$, and thus $\bp(W_\infty>0)=1$. The fact that
$p(\beta) = \frac{\beta^2 Q(0)}{2}$ is now easily seen from
(\ref{winfpbet}).

\hfill
$\Box$

\vspace{0.5cm}

Of course, Proposition \ref{pp} would be meaningless without some simple sufficient conditions on $Q$ ensuring hypothesis (H). Those sufficient conditions will be given in the following
\begin{proposition}
Assume $d\ge 3$, that $Q$ is a positive radial function from $\rd$ to $\R$, and write $Q(x)=\tq(|x|)$ for $x\in\rd$, where $\tq$ is a positive function from $\R$ to $\R$. Assume that $\beta$ is small enough and that
$$
\int_0^\infty x \tq(x) dx<\infty.
$$
Then hypothesis {\bf (H)} is satisfied.
\end{proposition}

\vspace{0.5cm}

\noindent {\bf Proof:}
We will recall first some results presented in \cite{Yo}:
let us denote by $R_d$ the Bessel process in dimension $d$, by
$\{ l_t^x(R_d); x>0, t\ge
0 \}$ the local time of the Bessel process, and by $X$ a standard
planar
Brownian Motion (we will assume that all those objects can be defined on
$(\hat\Omega,\hat\cf,\hat P)$). Then we have, for $d\ge 3$,
$$
\{ l_{\infty}^x (R_d); x >0\} \stackrel{(\cl)}{=}
\left\{
\frac{1}{(d-2)x^{d-3}}
\vert X_{x^{d-2}} \vert^2;  x >0 \right\}.
$$
Now, obviously, if $Q(x)=\tq(|x|)$, we have
$$
I_{\infty}(Q)=
 \int_0^\infty  \hat Q (R_d(s)) ds  = \int_0^\infty \hat Q(x)
 l_\infty^x (R_d) dx.
$$
Hence, by changing variables, we get
\begin{eqnarray*}
E_\omega  \left[ e^{\frac{\beta^2}{2} I_{\infty}(Q)} \right] & = &
\eo\lc \exp \lp \frac{\beta^2}{2}\int_0^\infty \hat Q(x)
 l_\infty^x (R_d) dx
\rp \rc\\
&=&\eo\lc \exp \lp
\int_0^\infty
\beta^2 S(v) |X_v|^2 dv
\rp \rc\\
&=&\prod_{i\le 2}
\eo\lc \exp \lp
\int_0^\infty
\beta^2 S(v) |X_v^i|^2 dv
\rp \rc,
 \end{eqnarray*}
where
$$
S(v)=\frac{\tq(v^{\frac{1}{d-2}})} {2(d-2)^2
v^{\frac{2(d-3)}{d-2}}}.
$$
Observe now that, following Fernique's definitions and results (see \cite{Fe}),
\begin{enumerate}
\item
$\vp\mapsto N(\vp)=(\int_0^\infty
S(v) \vp^2(v) dv)^{1/2}$ is a gauge on $C(\R_+)$ (see \cite[Definition 1.2.1]{Fe}). Indeed, the only fact that has to be checked is to show that $N$ is lower semi-continuous, ie that all the sets $\{\vp;N(\vp)\le M  \}$ are closed in $C(\R_+)$, for any value of $M\le 0$. But this point is a direct consequence of Fatou's lemma.
\item
Since, for $i=1,2$, $X_s^i$ is a Gaussian process, if $N(X^i)<\infty$ almost
 surely and $\beta$ is small enough, then $\eo[e^{ \beta^2 N(X^i)}]<\infty$  (cf \cite[Theorem 1.2.3]{Fe}).
\end{enumerate}
Thus, condition {\bf (H)} is now implied, for $\beta$ small enough, by the condition
$$
\hp\lp \int_0^\infty S(v) |X_v^1|^2 <\infty\rp=1.
$$
However, by \cite[Proposition 2.2]{PY}, this occurs iff $\int_0^\infty
v S(v)  dv<\infty$, which is equivalent to $\int_0^\infty x \tq(x) dx<\infty$
by an elementary change of variables.

\hfill
$\Box$

\vspace{0.5cm}

We can now state an improved concentration result
below $\frac{\beta^2Q(0)}{2}$ in the weak disorder regime:
\begin{proposition}\label{p15}
Assume (H) is satisfied, and that
$$
\eo\lc I_\infty(Q)\,e^{\frac{\beta^2}{2} I_\infty(Q)} \rc <
\infty.
$$
Then there exists a positive constant $K_1$ depending on $\beta$
and $Q$ such that
$$
\bp \lp \log(Z_t) \le \frac{\beta^2}{2} t Q(0) - u \rp \le K_1
\exp \lp - \frac{u^2}{K_1} \rp,
$$
for all $u,t>0$.
\end{proposition}

\vspace{0.5cm}

\noindent {\bf Proof:}
This proof will be again divided in two steps.

\vspace{0.5cm}

\noindent {\bf Step 1:} {\it Some moment inequalities.}

\noindent
Using (\ref{sqzt}) and (\ref{espht}) we have, under assumption (H),
$$
\frac{\be \left[ Z_t^2 \right ]}{\left( \be [Z_t] \right)^2 }
=E_\omega  \left[ e^{\beta^2 \int_0^t Q( \omega_s^1 - \omega_s^2)
ds} \right] \le K_1, $$ for some positive constant $K_1.$ Then
Paley-Zygmund's inequality gives us
$$
\sup_t \bp \left( Z_t \ge \frac12 \be[Z_t ] \right) \le \sup_t
\frac14 \frac{\left( \be [Z_t] \right)^2 }{\be \left( Z_t^2 \right
)} \ge K_2,$$
 for some positive constant $K_2.$
For $t> 0$, set $\itq=\iot Q(\omega_s^1-\omega_s^2 ) ds$.
Given another positive constant $K_3$, and recalling notation (\ref{avpol}), we are
 now able to compute
\begin{align*}
&  \bp \left( Z_t \ge \frac12 \be[Z_t ],  \langle \itq
 \rangle_t \le K_3  \right) \\
&   \ge \bp \left( Z_t \ge \frac12 \be[Z_t ],  E_{\omega }\left[
 \itq   e^{\beta
\sum_{i=1}^2
 \int_{0}^{t}B(ds,\omega^i_{s})} \right] \le \frac{K_3\left( \be[\zt] \right)^2 }{4}  \right) \\
&  \ge \bp \left( Z_t \ge \frac12 \be[Z_t ] \right) - 1 \\
& \quad + \bp \left( E_{\omega }\left[ \itq e^{\beta \sum_{i=1}^2
 \int_{0}^{t}B(ds,\omega^i_{s})} \right] \le \frac{K_3\left( \be[\zt] \right)^2}{4}
 \right).
\end{align*}
However, Chebychev's inequality yields
\begin{align*}
&  \bp \left( E_{\omega }\left[ \itq e^{\beta \sum_{i=1}^2
 \int_{0}^{t}B(ds,\omega^i_{s})} \right] > \frac{K_3\left( \be[\zt] \right)^2}{4}
 \right) \\
&   \le \frac{4}{K_3 \left( \be[\zt] \right)^2} \be \left[
E_{\omega }\left[  \itq
 e^{\beta \sum_{i=1}^2
 \int_{0}^{t}B(ds,\omega^i_{s})} \right] \right]\\
 &  = \frac{4}{K_3 \left( \be[\zt] \right)^2}
E_{\omega }\left[ \itq  \be \left[
e^{\beta \sum_{i=1}^2
 \int_{0}^{t}B(ds,\omega^i_{s})} \right] \right].
\end{align*}
Since
$$
\be \left[ e^{\beta \sum_{i=1}^2
 \int_{0}^{t}B(ds,\omega^i_{s})} \right] =
\left( \be[\zt] \right)^2 e^{\beta^2 \itq},
$$
we finally obtain
$$
\bp \left( Z_t \ge \frac12 \be[Z_t ],  \langle \itq \rangle_t \le
K_3  \right)\ge K_2 - \frac{4}{K_3} E_\omega \left[
\itq
e^{\beta^2 \itq} \right].
$$
So, our assumptions imply that, choosing $K_3$ large enough, we have
\begin{equation}\label{momineq}
\bp \left( Z_t \ge \frac12 \be[Z_t ],  \langle \itq \rangle_t \le K_3  \right) \ge
\frac{1}{K_3}.
\end{equation}

\vspace{0.5cm}

\noindent {\bf Step 2:} {\it Application of the concentration inequalities.}

\noindent
For a given Gaussian landscape $B$, that can be considered as an element of $M$, set
$$
Z_{t}(B)=E_{\omega }\left[ e^{ \beta
\iot B(ds,\om_s)}\right],
$$
and
$$
 \langle f(\omega^1,\omega^2) \rangle_t^B
 := \frac{  E_{\omega }\left( f(\omega^1,\omega^2 ) e^{\beta \sum_{i=1}^2
 \int_{0}^{t}B(ds,\omega^i_{s})} \right)}{Z_t^2(B)}.
 $$
For the constant $K_3$ used for inequality (\ref{momineq}), we can now consider the set
$$
A:= \lcl g\in M;\, Z_t(g) \ge \frac12 \be [Z_t],\, \langle \itq
\rangle_t^g \le K_3 \rcl,
$$
and we have
checked that
$$
\bp \left( B \in A \right) \ge \frac{1}{K_3}.
$$
Applying Lemma
\ref{ltala}, this yields that, for all $u>0$,
\begin{equation}\label{desbona}
\bp \left( q_A > u +K_4 \right) \le 2 \exp \left( - \frac{u^2}{2}
\right),
\end{equation}
with $K_4 = (2 \log(2K_3) )^\frac12.$

Consider now another Gaussian landscape $\bar B$, but keep the notation
$Z_t=Z_t(B)$. We can write
\begin{eqnarray}
Z_t & = & E_{\omega }\left[ e^{ \beta \int_{0}^{t} B(ds, \omega
_{s})}\right]  \nonumber\\
& = & E_{\omega } \left[ e^{ \beta \int_{0}^{t} B(ds, \omega
_{s})- \bar  B(ds, \omega _{s}) } e^{ \beta \int_{0}^{t} \bar
B(ds,
\omega _{s})}\right]  \nonumber\\
& = & Z_t( \bar B) \langle e^{ \beta \int_{0}^{t} B(ds, \omega
_{s})- \bar  B(ds, \omega _{s}) } \rangle_t^{\bar B} \nonumber\\
& \ge & Z_t( \bar B) e^{ \beta  \langle  \int_{0}^{t} B(ds, \omega
_{s})- \bar  B(ds, \omega _{s}) \rangle_t^{\bar B}},\label{minztztb}
\end{eqnarray}
where in the last step we have used Jensen's inequality. Suppose now that $B-\bar B=\tilde g$, where $\tilde g(t)=\iot g(s) ds$ and $g$ is the inverse Fourier transform of an element of $\ch$. Notice that $g$ admits the representation
$$
g(s,x)=\ird e^{\io\la x}h(s,\la) \hq(d\la), \, \mbox{ with }\,
|h|_{\ch}^{2}= \int_0^\infty \ird |h(s,\la)|^2 \hq(d\la)
ds<\infty.
$$
Furthermore,
\begin{multline}\label{difbbar}
\lln \lla \iot B(ds,\om_s)-\bar B(ds,\om_s)\rat^{\bar B} \rrn
=\lln  \iot\ird \lla e^{\io\la \om_s} \rat^{\bar B}h(s,\la) \hq(d\la) ds\rrn\\
\le |h|_{\ch}
\lp \iot\ird \lla e^{\io\la (\om_s^1-\om_s^2)} \rat^{\bar B}
\hq(d\la) ds\rp^{1/2}
= |h|_{\ch}
\lp \lla \itq \rat^{\bar B} \rp^{1/2}.
\end{multline}
Thus, putting together (\ref{minztztb}) and (\ref{difbbar}), we get, if $\bar B\in A$ and $B-\bar B=\tilde g$,
\begin{eqnarray*}
\log (Z_t)  & \ge & \log \lp Z_t( \bar B)\rp +  \beta \left \langle \int_{0}^{t}
B(ds, \omega _{s})- \bar  B(ds, \omega _{s}) \right\rangle_t^{\bar B} \\
& \ge & \log(\be[\zt]) - \log(2) - \beta\vert h \vert_\ch K_3^\frac12\\
&=&\frac{\beta^2 t Q(0)}{2}- \log(2) - \beta\vert h \vert_\ch K_3^\frac12.
\end{eqnarray*}
Obviously, one can choose, in the above inequality, the norm $\vert h \vert_\ch$ as close as desired to $q_A$. Thus, we get
$$
\log (Z_t)\ge \frac{\beta^2 t Q(0)}{2}- \log(2) - \beta q_A
K_3^\frac12,
$$
and using (\ref{desbona}) we have that, for all $u>0$,
the event
$$
\log (Z_t)\ge \frac{\beta^2 t Q(0)}{2}- \log(2) -  \beta K_3^\frac12 (u+K_4)
$$
holds with probability larger than $ 1 - \exp(-\frac{u^2}{2})$.
The proof is now easily completed.

\hfill $\Box$

\section{The strong disorder regime}

In this section, we will give some examples of Gaussian polymers in the strong disorder regime. We will begin with a general sufficient condition. Recall that $Q$ is the covariance of our noise $B$.
\begin{theorem}\label{t1}
Let $p>1$ be a
constant, $\{\Lambda_s ;s \in \R^+\}$ a family of subsets of $\rd$
and
$$
\kappa=\frac12 \beta^2 Q(0) (1-4q)^2 q^{-1},
$$ where $q$ is the conjugate exponent of $p$. Based on these notations, set
\begin{eqnarray*}
v(s)&=&\inf_{x \in \Lambda_s} Q(x)\\
w(s)&=&\left( \inf_{x \in \Lambda_s} Q(x) \right)
\hp^{\frac{1}{p}} (\omega_s^1 - \omega_s^2 \in \Lambda_s^c ) \, e^{\kappa s},
\end{eqnarray*}
and assume that
\begin{itemize}
\item[{\bf (H1)}]
$\int_0^\infty v(s)= \infty$ and $\int_0^\infty w(s)<\infty$.
\end{itemize}
Then
$$
\bp(W_\infty = 0)=1.
$$
\end{theorem}

\vspace{0.5cm}

\noindent {\bf Proof:} Since $W_\infty \ge 0$, we have, for any
$\theta>0$,
$$
\be[W_\infty^\theta] = \be[ \liminf_{t \to \infty} W_t^\theta ] \le
\liminf_{t \to \infty} \be[W_t^\theta ].
$$
Thus, it is enough to check
that
$$
\liminf_{t \to \infty} \be[W_t^\theta ]=0.
$$
Recall now the martingale decomposition we got for $W$ at (\ref{mgw}): setting
$$
X_s=\exp\lp N_s-\frac{\beta^2 Q(0)s}{2}\rp,
$$
one can write, for $t\ge 0$,
$$
\wt=
1+\beta
\iot\ird\eo
\lc  e^{\iota\la\om_s} X_s\rc
M(ds,d\la).
$$
Pick now $0<\theta<1$. An application of It\^o's formula gives
\begin{eqnarray*}
W_{t}^{\theta } &=&1+\beta \theta
\int_{0}^{t}\ird W_{s}^{\theta -1}E_{\omega }\left[ X_{s}e^{\io\la\om_s}
\right] M(ds,d\la) \\
&&-\frac{\beta ^{2}}{2}\theta (1-\theta
)\int_{0}^{t}\ird W_{s}^{\theta -2}\left( E_{\omega
}\left[ X_{s}e^{\io\la\om_s}\right] \right)^{2}\hq(d\la)ds.
\end{eqnarray*}
Then, taking expectations, we obtain
\begin{eqnarray*}
\be\left[ W_{t}^{\theta }\right] &=&1-\frac{\beta ^{2}}{2}\theta
(1-\theta
)\be\left[ \int_{0}^{t}\ird W_{s}^{\theta -2}\left( E_{\omega
}\left[ X_{s}e^{\io\la\om_s}\right] \right)^{2}\hq(d\la)ds\right] \\
&=&
1-\frac{\beta ^{2}}{2}\theta
(1-\theta
)\be\left[ \int_{0}^{t} W_{s}^{\theta -2}
\eo\lc X_s^1 X_s^2 Q(\om_s^1-\om_s^2) \rc
ds\right] .
\end{eqnarray*}
Hence
\begin{eqnarray}
&&\be\left[ W_{s}^{\theta -2}E_{\omega }\left[
X_{s}^{1}X_{s}^{2}Q(\omega
_{s}^{1}-\omega _{s}^{2})\right] \right]  \nonumber \\
&\geq &\left( \inf_{\Lambda _{s}}Q\right) \be\left[ W_{s}^{\theta
-2}E_{\omega }\left[ X_{s}^{1}X_{s}^{2}\1_{\left\{ \omega
_{s}^{1}-\omega
_{s}^{2}\in \Lambda _{s}\right\} }\right] \right]  \nonumber \\
&=&\left( \inf_{\Lambda _{s}}Q\right) \be \left[ W_{s}^{\theta
}\right]
-\left( \inf_{\Lambda _{s}}Q\right) \be \left[ W_{s}^{\theta -2}E_{\omega }%
\left[ X_{s}^{1}X_{s}^{2}\1_{\left\{ \omega _{s}^{1}-\omega
_{s}^{2}\in \Lambda _{s}^{c}\right\} }\right] \right] . \label{bb}
\end{eqnarray}%
On the other hand, H\"older's inequality yields, for any conjugate exponents $p,q$,
\begin{eqnarray}
&&\be \left[ W_{s}^{\theta -2}E_{\omega }\left[ X_{s}^{1}X_{s}^{2}\1_{\left%
\{ \omega _{s}^{1}-\omega _{s}^{2}\in \Lambda _{s}^{c}\right\} }\right] %
\right]  \nonumber \\
&=&E_{\omega }\left[ \1_{\left\{ \omega _{s}^{1}-\omega
_{s}^{2}\in \Lambda
_{s}^{c}\right\} }\be \left[ W_{s}^{\theta -2}X_{s}^{1}X_{s}^{2}\right] %
\right]  \nonumber \\
&\leq &\hp^{\frac{1}{p}}\left( \omega _{s}^{1}-\omega _{s}^{2}\in
\Lambda _{s}^{c}\right) E_{\omega }^{\frac{1}{q}}\left[
\be ^{q}\left[ W_{s}^{\theta -2}X_{s}^{1}X_{s}^{2}\right] \right]
.  \label{cc}
\end{eqnarray}
In particular, if $q=\theta^{-1}$, invoking the fact that
$E\left[ X^{\rho }\right] \leq E^{\rho }\left[ X\right] $
for $\rho \leq 1$ and $X\ge 0$, we get that
\begin{eqnarray}
&&E_{\omega }\left[ \be ^{q}\left[ W_{s}^{\theta -2}X_{s}^{1}X_{s}^{2}%
\right] \right]  \nonumber \\
&=&e^{-\frac{\beta ^{2}}{2}  sQ(0)}E_{\omega }\left[
\be ^{q}\left[ E_{\omega }^{\theta -2}\left[ e^{\beta
\int_{0}^{s}B(du,\omega _{u})}\right] e^{\beta
\int_{0}^{s}B(du,\omega _{u}^{1})}e^{\beta \int_{0}^{s}B(du,\omega
_{u}^{2})}\right] \right]  \nonumber \\
&\leq &e^{-\frac{\beta ^{2}}{2}  sQ(0)}E_{\omega }\left[
\be \left[
E_{\omega }^{q(\theta -2)}\left[ e^{\beta \int_{0}^{s}B(du,\omega _{u})}%
\right] e^{q\beta \int_{0}^{s}B(du,\omega _{u}^{1})+B(du,\omega _{u}^{2})}%
\right] \right] \label{dd} \\
&\le& e^{-\frac{\beta ^{2}}{2}  sQ(0)} \be \left[ \eo\lc e^{q\beta
(\int_{0}^{s}B(du,\omega _{u}^{1})+B(du,\omega _{u}^{2}))
+(1-2q)\beta \int_{0}^{s}B(du,\omega _{u}^{3})}\rc \right]\le e^{q
\kappa s}.\nonumber
\end{eqnarray}
Then, putting together (\ref{bb}), (\ref{cc}) and (\ref{dd}), we obtain
$$
\be \left[ W_{s}^{\theta -2}E_{\omega }\left[
X_{s}^{1}X_{s}^{2}Q(\omega
_{s}^{1}-\omega _{s}^{2})\right] \right]
\geq
v(s) \be \left[W_{s}^{\theta }\right] -w(s),
$$
Consequently,
\begin{eqnarray*}
\be \left[ W_{t}^{\theta }\right] &\leq &1-\frac{\beta
^{2}}{2}\theta (1-\theta )\int_{0}^{t} v(s)
 \be \left[
W_{s}^{\theta }\right] ds \\
&&+\frac{\beta ^{2}}{2}\theta (1-\theta )\int_{0}^{t}w(s) ds.
\end{eqnarray*}
However, using our assuptions {\bf (H1)} and  setting
$$
\gamma=\frac{\beta ^{2}}{2}\theta (1-\theta ),
\quad \mbox{ and }\quad
\delta=1+\gamma\int_0^{\infty} w(s) ds,
$$
we get
$$
\be \left[ W_{t}^{\theta }\right] \le \delta- \gamma \iot  v(s)
\be \left[ W_{s}^{\theta }\right] ds,
$$
and by a standard comparison argument for ordinary differential equations, this yields
$$
\be \left[ W_{t}^{\theta }\right] \le \delta e^{-\gamma \iot v(s)
ds},
$$
and hence, invoking again Hypotesis (H1),
$$
\lim_{t\to\infty}\be \left[ W_{t}^{\theta }\right] =0,
$$
which proves our claim.

\hfill $\Box$

\medskip
\begin{example} Consider $d \ge 1$, and assume that the covariance function
$Q$ satisfies
\begin{equation}\label{exemple}
c_1 (1 + |x|^2)^{- \lambda} \le Q(x) \le c_2 (1+ |x|^2)^{-\hat
\lambda}, \end{equation}
 for some constants
 $c_1>0, c_2 >0 $ and $0 < \hat\lambda \le \lambda < \frac12$. Then the polymer
 will be in the strong disorder
 regime for any value of $\beta>0$.
\end{example}

\vspace{0.5cm}

\noindent {\bf Proof:} Observe that there exist some positive
definite functions $Q$ satisfying (\ref{exemple}), since a
function of the type $c_1 (1 + |x|^2)^{- \lambda}$ is the Fourier
transform of a tempered measure (see \cite[page 288]{Gel}).

Now, Theorem \ref{t1} can be applied with an arbitrary constant
$p>1$, by choosing the set $\Lambda_s$ as the centered ball of
radius $s^\alpha$ in $\R^d$, with $\al >1$ such that $\al \lambda
< \frac12.$ Indeed, it is easily seen in this case that
$$
v(s) \ge   c_1 ( s^{2 \alpha} + 1)^{-\lambda},
$$
and
$$w(s) \le c_2 ( s^{2 \alpha} + 1)^{- \hat \lambda} \exp \left( - \frac{s^{2 \alpha-1}}{8p} +
\kappa s \right),
$$
which proves that the assumption {\bf (H1)} is verified.

\hfill $\Box$


\begin{thebibliography}{99}
\bibitem{AZ} Albeverio, S; Zhou, X. (1996). A martingale approach to directed polymers in a random environment.  {\it J. Theoret. Probab.}  {\bf 9}  ,  no. 1, 171--189.

\bibitem{Bo}  Bolthausen, E. (1989). A note on the diffusion of directed polymers in a random environment.  {\it Comm. Math. Phys.}  {\bf 123}  (1989),  no. 4, 529--534.


\bibitem{CH}  Carmona, P; Hu, Y. (2002). On the partition function of
a directed polymer in a Gaussian random environment.
\textit{Probab. Theory Relat. Fields,} \textbf{124}, 431-457.

\bibitem{CM} Carmona, R; Molchanov, S. A. (1994). Parabolic Anderson problem and intermittency.  {\it Mem. Amer. Math. Soc.}  {\bf 108}.

\bibitem{CV} Carmona,R; Viens, F. (1998). Almost-sure exponential
behavior of a stochastic Anderson model with continuous space
parameter. \textit{Stochastics and Stochastic Reports,}
\textbf{62}, 251-273.

\bibitem{CY}  Comets, F; Yoshida, N. (2003). Brownian Directed Polymers
in Random Environment. \textit{Preprint}.

\bibitem{CO}Conlon, J; Olsen, P. (1996). A Brownian motion version of the directed polymer problem.  {\it J. Statist. Phys.}  {\bf 84}  ,  no. 3-4, 415--454.

\bibitem{Co}Coyle, L. (1996). A continuous time version of random walks in a random potential.  {\it Stochastic Process. Appl.}  {\bf 64}  ,  no. 2, 209--235.

\bibitem{CS} Cranston, M; Mountford, T; Shiga, T. (2002). Lyapunov exponents for the parabolic Anderson model.  {\it Acta Math. Univ. Comenian. (N.S.)}
{\bf 71},  no. 2, 163--188.

\bibitem{DS}Dawson, D.A.; Salehi, H. (1980).
Spatially homogeneous random evolutions. {\it J. Multivariate
Anal.} {\bf 10}, no. 2, 141-180.

\bibitem{DSp}Derrida, B; Spohn, H. (1988). Polymers on disordered trees, spin glasses, and traveling waves. New directions in statistical mechanics (Santa Barbara, CA, 1987).  {\it J. Statist. Phys.}  {\bf 51}  ,  no. 5-6, 817--840. Derrida Spohn


\bibitem{Fe}Fernique, X. (1997). Fonctions al\'eatoires
gaussiennes, vecteurs al\'eatoires gaussiens. Centre de Recherches
Math\'ematiques, Montreal.

\bibitem{FU}  Feyel, D; \"Ust\"unel, A.S. (2004). Monge-Kantorovich Measure
Transportation and Monge-Amp\`ere Equation on Wiener Space.
\textit{Probab. Theory Relat. Fields,} \textbf{128}, 347-385.

\bibitem{Gel} Gel'fand, I.M.; Vilenkin, N. Ya.(1964).
{\it Generalized Functions. Volume I}. Academic Press, New York.

\bibitem{IS}Imbrie, J; Spencer, T. (1988).
Diffusion of directed polymers in a random environment.
{\it J. Statist. Phys.} {\bf 52}, no. 3-4, 609--626.

\bibitem{Ma}Major, P. (1981). {\it Multiple Wiener-Itô integrals. With applications to limit theorems.} Lecture Notes in Mathematics, {\bf 849}.
Springer, Berlin.

\bibitem{Mal} Malliavin, P. (1997) {\it Stochastic Analysis}.
Springer, Berlin.

\bibitem{Nu} Nualart, D. (1995) {\it The Malliavin Calculus and
Related Topics}. Springer, Berlin.

\bibitem{PY}  Pitman, J; Yor, M.  (1982). A decomposition of Bessel
bridges. \textit{Z. Wahrscheinlichkeitstheorie verw. Gebiete}
\textbf{59}, 425-457.

\bibitem{RY}Revuz,D;  Yor, M. (1991). {\it Continuous martingales and Brownian
motion}. Springer, Berlin.

\bibitem{SY}  Salminen, P; Yor, M. (2003). Properties of perpetual integral
functionals of Brownian motion with drift. \textit{Preprint}.

\bibitem{Sn}Sznitman, A. (1998). {\it Brownian motion, obstacles and random media}. Springer Monographs in Mathematics. Springer-Verlag, Berlin.

\bibitem{Si}Sinai, Yakov G. (1995). A remark concerning random walks with random potentials.  {\it Fund. Math.}  {\bf 147},  no. 2, 173--180.


\bibitem{Tbk} Talagrand, M. (2003). {\it Spin Glasses: A challenge for
Mathematicians.} Springer, Berlin.

\bibitem{TV1} Tindel, S; Viens, F. (2002). Almost sure exponential behaviour for a parabolic SPDE on a manifold.  {\it Stochastic Process. Appl.}  {\bf 100}  , 53--74.

\bibitem{TV2}Tindel, S; Viens, F. (2001). Relating the almost-sure Lyapunov exponent of a parabolic SPDE and its coefficients' spatial regularity. To appear at {\it Potential Anal.}

\bibitem{U}\"Ust\"unel, A.S. (1995). {\it An introduction to analysis on Wiener space}.
 Lecture Notes in Mathematics {\bf 1610}. Springer-Verlag, Berlin.

\bibitem{U2} \"Ust\"unel, A.S. {\it Private Communication}.

\bibitem{Wu}W\"uthrich, M. (1998). Superdiffusive behavior of two-dimensional Brownian motion in a Poissonian potential.  {\it Ann. Probab.}  {\bf 26}  ,  no. 3, 1000--1015.

\bibitem{Yo} Yor, M. (1992). {\it Some aspects of Brownian motion. Part
I. Some special functionals.} Lectures in Mathematics ETH
Z\"urich. Birk\"auser, Basel.

\end{thebibliography}
\end{document}